\documentclass[11pt,a4paper,leqno]{article}
\usepackage{amsmath}
\usepackage{amssymb}
\usepackage{color}
\usepackage{yfonts}
\advance\topmargin by -2.2cm
\newcommand{\bbr}{I\!\!R}

\newcommand{\bbn}{I\!\!N}

\newcommand{\bbf}{I\!\!F}

\newcommand{\calb}{{\cal B}}

\newcommand{\cald}{{\cal D}}
\newcommand{\cale}{{\cal E}}
\newcommand{\calf}{{\cal F}}
\newcommand{\calg}{{\cal G}}

\newcommand{\call}{{\cal L}}

\newcommand{\caln}{{\cal N}}

\newcommand{\barr}{\begin{array}}
\newcommand{\earr}{\end{array}}
\newcommand{\beqq}{\begin{equation}}
\newcommand{\eeqq}{\end{equation}}
\newcommand{\beao}{\begin{eqnarray*}}
\newcommand{\eeao}{\end{eqnarray*}\noindent}
\newcommand{\beam}{\begin{eqnarray}}
\newcommand{\eeam}{\end{eqnarray}\noindent}

\newcommand{\halmos}{\quad\hfill\mbox{$\Box$}}

\newcommand{\la}{\lambda}

\newcommand{\si}{\sigma}
\newcommand{\al}{\alpha}
\newcommand{\vth}{\vartheta}

\newcommand{\om}{\omega}

\newcommand{\vep}{\varepsilon}

\newcommand{\wh}{\widehat}
\newcommand{\wt}{\widetilde}

\newcommand{\lra}{\longrightarrow}

\newcommand{\nto}{n\to\infty}

\begin{document}
%%%%%%%%%%%%%%%%%%%%%%%%%%%%%

{\huge\bf LAMN in a class of parametric models \\for null recurrent diffusions} \\  

{\bf November 6, 2017
%\today
}\\

{\bf Reinhard H\"opfner \\Carina Zeller \\ Johannes Gutenberg Universit\"at Mainz}\\

\vskip1.5cm\small
{\bf Abstract: } We study statistical models for one-dimensional diffusions which are recurrent null. A first parameter in the drift is the principal one, and determines regular varying rates of convergence for the score and the information process. A finite number of other parameters, of secondary importance, introduces additional flexibility for the modelization of the drift, and does not perturb the null recurrent behaviour. Under time-continuous
 observation we obtain local asymptotic mixed normality (LAMN), state a local asymptotic minimax bound, and specify asymptotically optimal estimators.    
 
\\
{\bf AMS Mathematics Subjects Classifications: } 62 F 12, 62 M 05, 60 J 60  \\
\normalsize

\vskip1.2cm
%%%%%%%%%%%%%%%%%%%%%%%%%%%%%%%%%%%%%%%%%%%%%%%%
%  introduction 
%%%%%%%%%%%%%%%%%%%%%%%%%%%%%%%%%%%%%%%%%%%%%%%%%%

In a joint paper with Kutoyants \cite{HK 03} we considered a one-parametric statistical model 
\beqq\label{HK_model_1}
d\xi_t \;=\; \vth\, f_1(\xi_t)\, dt \;+\; \si\, dW_t  \quad,\quad t\in[0,n] 
\eeqq
where the path of the diffusion $\xi$ is observed continuously in time, and where the shape of the function $f_1$ in the drift combined with the range of the parameter $\vth$
\beqq\label{2_intro}
f_1(x) \;=\; \frac{x}{1+x^2} \quad,\quad \Theta_1 := (-\frac12\si^2 ,+\frac12\si^2) \quad,\quad  \vth\in \Theta_1   
\eeqq
makes the diffusion $\xi=(\xi_t)_{t\ge 0}$ in \eqref{HK_model_1} recurrent null. Relying on limit theorems for additive functionals and for martingale additive functionals from H\"opfner and L\"ocherbach \cite{HL 03}, we proved LAMN and specified optimal  estimators in \cite{HK 03}, with speed of convergence  
\beqq\label{speed_conv}
\mbox{ $\sqrt{ n^{\al(\vth)} }\;\;$ with $\;\;\al(\vth):=\frac12(1-\frac{2\vth}{\si^2}) \;$}   
\eeqq
such that $\vth\to\al(\vth)$ maps $\Theta_1$ onto $(0,1)$. The diffusion model \eqref{HK_model_1} with $f_1$ and $\Theta_1$ specified by \eqref{2_intro} is a prototype example for classes of statistical models where at every point of the model, both the speed of convergence of optimal estimators and the structure of their limit distributions are different (similiar examples for Markov step processes in  \cite{Ho 90}, \cite{Ho 93}%, example 2.23 b)
). The speed of convergence in \eqref{speed_conv} corresponds to the index of regular variation of the tails of suitably defined 'life cycles' of the process, i.e.\ suitably defined independent excursions away from state $0$.

\cite{HK 03} also considered the case of drift contamination by unknown smooth perturbations $g$ ranging over some class $\calg_C:=\{ g:\bbr\to\bbr \;\,\mbox{smooth such that}\;\, {\rm supp}(g)\subset C \}$. It was assumed that the compact set $C$ supporting the perturbations is fixed and known, so the model is   
\beqq\label{HK_model_2}
d\xi_t\; =\; \left( \vth f_1 + g \right)(\xi_t)\, dt \;+\; \si\, dW_t  \quad,\quad t\in[0,n]  \quad,\quad   \vth\in\Theta_1 \quad,\quad g\in\calg_C \;. 
\eeqq
Thanks to $\,{\rm supp}(g)\subset C\,$ we still can estimate $\,\vth\,$ from observation of $\xi$ restricted to the complement $C^c$. Such estimators conserve the optimal speed of convergence \eqref{speed_conv} from model \eqref{HK_model_1} even if limit distributions are more spread out.  If we focus now on one particular  $g\in\calg_C$ which we assume fixed and known, we can take such estimators as preliminary estimators to be modified by LeCam's one step correction in LAMN setting, and thus arrive at explicit optimal estimators for the unknown parameter.

Beyond \eqref{HK_model_1} and \eqref{HK_model_2}, the present paper aims at $(1+m)$-dimensional parametric models for null recurrent diffusions of the following type. Keeping $f_1$ and $\Theta_1$ as above in \eqref{2_intro}, we consider     
\beqq\label{1}
d\xi_t = \left( \vth_1 f_1 + \sum_{\nu=1}^m \vth_{2,\nu} f_{2,\nu} \right)\!(\xi_t)\, dt + \si dW_t  \quad,\quad t\in[0,n]  \quad,\quad   \vth_1\in\Theta_1  \quad,\quad  \vth_{2,\nu}{\in}\bbr 
\eeqq
with a finite number of functions  $\,f_{2,\nu}:\bbr\to\bbr\,$, fixed and known and with suitable properties,  whose presence in the drift does not perturb the null recurrent behaviour of $\xi=(\xi_t)_{t\ge 0}$. We are not interested in functions $f_{2,\nu}$ having compact support. The task is to estimate the unknown parameter 
$$
\vth := (\vth_1,\vth_{2,1},\ldots,\vth_{2,m}) \;\in\; \Theta := \Theta_1\times\bbr^m
$$
from continuous-time observation of $\xi$ over a long time interval.

{\bf Example 1: } For $m=1$, we might consider $f_2(x):=\frac{\sin x}{x}$ in \eqref{1}, $x\in\bbr$. 
\\
For $m=2\ell$ with arbitrary $\ell\in\bbn$, we might define functions $f_{2,\nu}$ by 
\beqq\label{example}
f_{2,2k}(x) \;:=\; f_1(x)\cos (k x) \quad,\quad f_{2,2k+1}(x):=f_1(x)\sin (k x) \quad,\quad  1\le k\le \ell \;,  
\eeqq
and thus consider in the drift of \eqref{1} finite Fourier series of order~$\ell\,$ tempered by $f_1$. 

\vskip0.5cm
This example indicates that the type of models \eqref{1} should be flexible enough to take into account fine differences in nullrecurrent behaviour over a broad range \eqref{speed_conv} of possible speeds of convergence.

For such models, we shall prove LAMN in the $(1{+}m)$-dimensional parameter $\vth\in\Theta$, and specify optimal estimators in the sense of the local asymptotic minimax theorem. The main results are lemma 1 (null recurrence and invariant density) in section \ref{assumptions}, theorem 1 (weak convergence of the score function martingale together with the information process) in section \ref{convergence theorem}, and theorem 2 (local asymptotic minimax bound) in section \ref{comp_est}.

The notions of local asymptotic normality (LAN) or local asymptotic mixed normality (LAMN) go back to LeCam, see LeCam \cite{LC 69}, H\'ajek \cite{Ha 70}, Jeganathan \cite{Je 82}, Davies \cite{Da 85}, Jeganathan \cite{Je 88}, and LeCam and Yang \cite{LY 90}. Ibragimov and Khasminskii \cite{IH 81} (see p.\ 162 there) and Kutoyants \cite{Ku 04} (see p.\ 114 there) take a different route towards a local asymptotic minimax theorem under LAN: their bound uses neighbourhoods of the true parameter which are small but not shrinking at a rate which corresponds to the norming sequence for the score under the true parameter. For results and proofs in the LAMN context we refer to sections 5.1, 6.1 and 7 of \cite{Ho 14}.

%%%%%%%%%%%%%%%%%%%%%%%%%%%%%%%%%%%%%%%%%%%%%%%%%%
%\newpage
%\vskip0.8cm
\section{The setting}\label{setting}
%%%%%%%%%%%%%%%%%%%%%%%%%%%%%%%%%%%%%%%%%%%%%%%%%%%

We discuss the probabilistic background for the diffusion \eqref{1} and specify the assumptions which will in force throughout the paper. Assumptions and results are stated in a first subsection, proofs and additional remarks in a second one.

%%%%%%%%%%%%%%%%%%%%%%%%%5
\subsection{Assumptions, null recurrence, likelihoods, some estimators}\label{assumptions}
%%%%%%%%%%%%%%%%%%%%%%%%%%%%

Throughout this paper, we consider the model \eqref{1}, with $f_1$ and the associated parameter given by 
\beqq\label{2}
f_1(x) \;=\; \frac{x}{1+x^2} \quad,\quad  \vth_1\in \Theta_1\quad,\quad \Theta_1 := (-\frac12\si^2 ,+\frac12\si^2)   \;. 
\eeqq
The functions $f_{2,1},\ldots,f_{2,m}$ in equation \eqref{1} satisfy the following.

{\bf Assumption 1: } The functions $f_{2,\nu}:\bbr\to\bbr$ in \eqref{1} are Lipschitz continuous and such that\\  
i) finite limits do exist for  $\,F_{2,\nu}(x):=\int_0^x f_{2,\nu}(y) dy\,$ as $x\to\pm\infty$, and are denoted by $\,F_{2,\nu}(\pm\infty)\in\bbr $; \\
ii)  for $\vep>0$ arbitrarily small, $\;\int_{\{|x|>1\}} [f_{2,\nu}]^2(x)\, |x|^{1-\vep}\, dx\;<\, \infty $; \\
iii) the functions occurring in \eqref{1} 
%$f_1, f_{2,1}, \ldots, f_{2,m}$ 
are linearly independent in the following sense: \\for $U$ open in $\bbr$ and real constants  $\zeta_1,\zeta_{2,1},\ldots,\zeta_{2,m}$,  
any representation $\;0 \,\equiv\, \zeta_1 f_1 + \sum\limits_{\nu=1}^m \zeta_{2,\nu} f_{2,\nu}\;$ on $U$ implies $\,0\,=\,\zeta_1{=}\zeta_{2,1}{=}\ldots{=}\zeta_{2,m}\,$.

%\vskip0.5cm
Note that the functions $\,f_2\,$ or $\,f_{2,1},\ldots,f_{2,2\ell}\,$ considered in example 1 above do satisfy all requirements of assumption~1:  iii) holds in both cases $m=1$ or $m=2\ell$, ii) is obvious, and elemen\-tary arguments (such as  $\,|\int_s^t \frac{\sin y}{y}\, dy\, | < \frac2s\,$ for all $0<s<t<\infty$ from \cite{He 93} Ch.~87) allow to check i). We put 
\beqq\label{3}
\Theta \;:=\; \Theta_1\times\bbr^m \;\;,\;\; \vth=(\vth_1,\vth_{2,1},\ldots,\vth_{2,m}) \quad,\quad  
\la_1(\vth) \;:=\; \frac{2}{\si^2}\vth_1 \;\;,\;\;
\la_{2,\nu}(\vth) \;:=\; \frac{2}{\si^2}\vth_{2,\nu} \;. 
\eeqq

{\bf Lemma 1: } Under assumption 1, the diffusion \eqref{1} is recurrent in the sense of Harris for every $\vth\in\Theta$.  With notation  \eqref{3}, the invariant measure of $\xi=(\xi_t)_{t\ge 0}$, unique up to constant multiples, is given by
$$
\mu^\vth(dx) \;=\; \frac{1}{\si^2}\, (\sqrt{1+x^2}\,)^{\la_1(\vth)}\, \exp\left\{\, \sum_{\nu=1}^m \la_{2,\nu}(\vth)\, F_{2,\nu}(x) \right\}\; dx    \quad,\quad \vth\in\Theta 
$$
on $(\bbr,\calb(\bbr))$. Here $\mu^\vth(\bbr)=\infty$, thus null recurrence holds for every $\vth\in\Theta$. \\

We strengthen that the Lebesgue density of the invariant measure $\mu^\vth$ varies regularly as $x\to+\infty$ and as $x\to-\infty$. The index of regular variation $\la_1(\vth)$, the same on the left and the right branch, depends on $\vth_1$ only and ranges over the interval $(-1,1)$. We have no regular variation for the functions $f_{2,\nu}$ in the drift, cf.\ example \eqref{example}, whereas they contribute to asymptotic constants in virtue of assumption~1~i). We have choosen $\Theta_1$ as the maximal open interval on which null recurrence holds.

Throughout this paper, the starting point $x_0\in\bbr$ for the diffusion \eqref{1} does not depend on $\vth\in\Theta$ and will be fixed. From now on, these assumptions and notations will remain in force (and we omit to recall this in the results below).

For the theoretical background of what follows we refer to the classical books  Liptser and Shiryaev \cite{LS 78} and Jacod and Shiryaev \cite{JS 87}; see also Ibragimov and Khasminskii \cite{IH 81}, Kutoyants \cite{Ku 04}, or section 6.2 in H\"opfner \cite{Ho 14}. Let $Q^\vth$ denote the law of  $\xi=(\xi_t)_{t\ge 0}$ under $\vth$ on the canonical path space $(D,\cald, \bbf)$ for cadlag processes (\cite{JS 87} Chapter VI), with canonical filtration 
$$
\bbf=(\calf_t)_{t\ge 0} \quad,\quad \calf_t \;=\; \bigcap\limits_{r>t} \si( \eta_s : 0\le s\le r ) \;. 
$$
Write $\eta=(\eta_t)_{t\ge 0}$ for the canonical process, i.e.\ the process of coordinate pro\-jec\-tions $\eta_t(f)=f(t)$, $f\in D$, and $\,m^{(\vth)}\,$ for the $Q^\vth$-martingale part of $\,\eta\,$. We introduce a $(1+m)$-dimensional $Q^\vth$-martingale $\,S(\vth) = \left( S(\vth)_t \right)_{t\ge 0}\,$:
\beqq\label{score}
S(\vth)_t = \left(\begin{array}{l} S^{(1)}(\vth)_t \\ S^{(2,1)}(\vth)_t \\ \quad\quad\vdots \\ S^{(2,m)}(\vth)_t \end{array}\right) \;\;:=\;\;  
\left(\begin{array}{l} 
\int_0^t \frac{1}{\si^2}\, f_1(\eta_s)\, dm^{(\vth)}_s  \\ 
\int_0^t \frac{1}{\si^2}\, f_{2,1}(\eta_s)\, dm^{(\vth)}_s \\
\quad\quad\quad\vdots \\
\int_0^t \frac{1}{\si^2}\, f_{2,m}(\eta_s)\, dm^{(\vth)}_s
\end{array}\right) \;,\quad t\ge 0   
\eeqq
where $m^{(\vth)}$ denotes the $Q^\vth$-martingale part of $\,\eta\,$. $\,S(\vth)$ will be called score martingale. The angle bracket $\left\langle\, S(\vth) \,\right\rangle$ of $S(\vth)$ under $Q^\vth$ is the process $J=(J_t)_{t\ge 0}$  
\beqq\label{inf_proc}
J_t = 
\left(\begin{array}{llll} 
\int_0^t \frac{1}{\si^2}\, [f_1]^2(\eta_s)\, ds  &  \ldots & \int_0^t \frac{1}{\si^2}\, [f_1 f_{2,\nu'}](\eta_s)\, ds & \ldots \\
\quad\quad\vdots  &  \ddots  &    &       \\
\int_0^t \frac{1}{\si^2}\, [f_{2,\nu} f_1](\eta_s)\, ds    &   &  \int_0^t \frac{1}{\si^2}\, [ f_{2,\nu} f_{2,\nu'} ](\eta_s)\, ds  &   \\
\quad\quad\vdots  &   &     &  \ddots     
\end{array}\right)  \;,\quad t\ge 0 
\eeqq
(with $\nu,\nu'$ ranging over $1,\ldots,m$) which is observable; $J$ will be called information process.  \\

{\bf Lemma 2: } For all pairs $\vth', \vth$ in $\Theta$, the laws $Q^{\vth'}, Q^\vth$ are locally equivalent relative to $\bbf$. The log-likelihood ratio process of $Q^{\vth'}$ with respect to $Q^\vth$ relative to $\bbf$ is given by 
\beqq\label{log-LR}
\Lambda^{ \vth' / \vth }_t \quad:=\quad 
(\vth'-\vth)^\top S(\vth)_t \;-\; \frac12 (\vth'-\vth)^\top J_t\; (\vth'-\vth)  \;,\quad t\ge 0 
\eeqq
where $\,\cdot\cdot^\top$ is the scalar product in $\bbr^{1+m}$. The information  $\,J_t\,$,  $0<t<\infty$, takes values in the set ${\tt D}^+$ of all strictly positive definite symmetric $(1{+}m){\times}(1{+}m)$ matrices, $Q^\vth$-almost surely for every $\vth\in\Theta$. \\

There are many possibilities to define estimators for the unknown parameter $\vth\in\Theta$ based on time-continuous observation of the diffusion path \eqref{1} on $[0,t]$.

{\bf Proposition 1: } We have well defined maximum likelihood (ML) estimators, of form 
$$
\wh \vth _t  \;=\;   \left(\begin{array}{l} {\wh \vth}^{(1)}_t \\ {\wh \vth}^{(2,1)}_t \\ \quad\vdots \\ {\wh \vth}^{(2,m)}_t \end{array}\right)  
\;:=\;  1_{ \{\,  J_t \,\in\, {\tt D}^+ \} }\,   \left[\, J_t \,\right]^{-1}\; Y_t \quad,\quad  
Y_t \;=\; \left(\begin{array}{l} Y^{(1)}_t \\ Y^{(2,1)}_t \\ \quad\vdots \\ Y^{(2,m)}_t \end{array}\right) \quad,\quad 0<t<\infty
$$
%for the unknown parameter $\vth\in\Theta$, 
%based on observation of the diffusion path \eqref{1} up to time $t$, 
where $Y^{(1)}$, $Y^{(2,\nu)}$ denotes common determinations  
\beao
& (t,\om) \to Y^{(1)}(t,\om) & \mbox{ for the family of stochastic integrals $\int  \frac{1}{\si^2}\, f_1(\eta_s)\, d\eta_s$ under $Q^\vth$ } \\
& (t,\om) \to Y^{(2,\nu)}(t,\om) & \mbox{ for the family of stochastic integrals $\int_0^t  \frac{1}{\si^2}\, f_{2,\nu}(\eta_s)\, d\eta_s $ under $Q^\vth$} 
\eeao
valid jointly under all $\vth\in\Theta$. For every $\vth\in\Theta$, the ML estimation error is such that  
\beqq\label{estim_errors_1}
\wh \vth _t - \vth \;\;=\;\; 1_{ \{  \,J_t \,\in\, {\tt D}^+ \} }\, [\, J_t \,]^{-1}\, S(\vth)_t \quad\quad\mbox{$Q^\vth$-almost surely}\;. 
\eeqq

\vskip0.5cm 
One might prefer to restrict observation of a null recurrent process to some fixed (and sufficiently large) compact set $A$ in $\bbr$, and thus consider estimators of the following type.

{\bf Proposition 2: } Fix $A$ compact in $\bbr$ such that the starting point $x_0$ of the diffusion is an interior point of $A$. Define $S(\vth,A)$, $J(A)$, $Y(A)$  in analogy to \eqref{score}, \eqref{inf_proc} and $Y$ above, but with $f_1$, $f_{2,1},\ldots,f_{2,m}$ replaced by $f_1 \cdot 1_A$, $f_{2,1} \cdot 1_A,\ldots, f_{2,m} \cdot 1_A$. Then $\,J(A)_t$ for $0<t<\infty$ takes values in ${\tt D}^+$, 
$\,Q^\vth$-almost surely for every $\vth\in\Theta$, and the estimator  
$$
\wt \vth (A)_t \;:=\;  1_{ \{  J(A)_t \,\in\, {\tt D}^+ \} }\,   \left[\, J(A)_t \,\right]^{-1}\;  Y(A)_t
$$
admits a representation 
\beqq\label{estim_errors_2}
\wt \vth (A) _t - \vth \;\;=\;\; 1_{ \{  J(A)_t \,\in\, {\tt D}^+ \} }\, [J(A)_t]^{-1}\, S(\vth,A)_t \quad\quad\mbox{$Q^\vth$-almost surely}
\eeqq
for every $\vth\in\Theta$. Under $\vth\in\Theta$, this estimator replaces the log-likelihood surface $\vth' \to \Lambda^{ \vth' / \vth }_t$ considered in \eqref{log-LR} by a modified (again 'inverse bowl-shaped') surface 
$$ 
\Theta \;\ni\; \vth' \quad\lra\quad
(\vth'-\vth)^\top S(\vth,A)_t \;-\; \frac12 (\vth'-\vth)^\top J(A)_t\; (\vth'-\vth)  \;. 
$$

\vskip0.5cm 
For the principal parameter $\vth_1\in\Theta_1$ we might think of an estimator which was optimal in the one-dimensional model \eqref{HK_model_1} considered in section 2 of \cite{HK 03}.

{\bf Remark 1: } Write $\,\check J _t = \int_0^t \frac{1}{\si^2}[f_1]^2(\eta_s) ds\,$ for the $(1,1)$-entry of $J_t$ in \eqref{inf_proc}, $\,\check S(\vth) _t$ for the first component of $S(\vth) _t$ in \eqref{score}, and $\check Y _t$ for the first component of $Y _t$ defined above. Consider   
$$
\check \vth (t) \;\;:=\;\; [\,\check J _t\,]^{-1}\,  \check Y _t \;\;=\;\;  \frac{ \int_0^t f_1(\eta_s) d\eta_s }{ \int_0^t [f_1]^2(\eta_s) ds } \quad,\quad 0<t<\infty    
$$
as an estimator for $\vth_1\in\Theta_1$. In our model \eqref{1}, this estimator is inconsistent: it admits a represention
$$
\check \vth (t) - \vth_1 \;\;=\;\; \check b _t (\vth) \;+\; [\,\check J _t\,]^{-1}  \check S(\vth) _t 
$$
under $Q^\vth$, with additional terms $\check b _t (\vth)$ which by the ratio limit theorem    
$$
\check b _t (\vth) \;\;=\;\;  \sum_{\nu=1}^m \vth_{2,\nu}  \frac{ \int_0^t \,[ f_1 f_{2,\nu} ](\eta_s)\, ds }{ \int_0^t [f_1]^2(\eta_s) ds } 
\quad\lra\quad 
\sum_{\nu=1}^m \vth_{2,\nu} \frac{ \mu^{\vth}\left( f_1  f_{2,\nu} \right) }{ \mu^{\vth}\left( f_1^2 \right) } \quad\mbox{ as $t\to\infty$}
$$  
converge $Q^\vth$-almost surely as $t\to\infty$.

%%%%%%%%%%%%%%
\subsection{The proofs}\label{stetting_proofs}
%%%%%%%%%%%%%

{\bf Proof of lemma 1: } The proof uses classical arguments, see e.g.\ Gihman and Skorohod \cite{GS 72} or Khasminskii \cite{Kh 80}; a resum\'e is in \cite{Ho 14} sections 9.2--9.4. Fix $\vth$, write for short $b$ for the drift in equation \eqref{1} under $\vth$
$$
b(x) \;=\; \left( \vth_1 f_1 + \sum_{\nu=1}^m \vth_{2,\nu} f_{2,\nu} \right)\!(x) \quad,\quad x\in\bbr \;, 
$$
and define functions $s$, $S$ on $\bbr$ by 
\beqq\label{harmonic}
S(x) := \int_0^x s(y)\, dy  \quad,\quad   s(y) := \exp\left( -\int_0^y \frac{2\, b}{\si^2}(v)\, dv \right) \;. 
\eeqq
%$S$ is harmonic for the Markov generator of the diffusion \eqref{1}: on $(D,\cald)$, $\,\wt\eta:=\left(S(\eta_t)\right)_{t\ge 0}$ is a local $(Q_\vth,\bbf)$-martingale; $\,\wt\eta\,$  solves the equation $\,d\wt\eta_s = \wt\si(\wt\eta_s)ds + dW_s\,$ with $\,\wt\si(y) = \si\cdot s (S^{-1}(y))\,$ where $S^{-1}$ denotes the function inverse of $S$.  
\eqref{2} and \eqref{3} together yield
\beqq\label{klein_s}
s(y) \;=\; \left( \sqrt{1+y^2}\, \right) ^{-\la_1(\vth)} \exp\left( -\sum_{\nu=1}^m \la_{2,\nu}(\vth) F_{2,\nu}(y) \right) \quad,\quad y\in\bbr
\eeqq
where $-1 < \la_1(\vth) < 1$, and we deduce from assumption 1 that 
\beqq\label{gross_S}
S(x) \;\sim\; \mbox{sgn}(x) \; \; \frac{1}{1-\la_1(\vth)}\, |x|^{1-\la_1(\vth)}\; \exp\left( -\sum_{\nu=1}^m \la_{2,\nu}(\vth) F_{2,\nu}(\pm\infty) \right)  \quad\mbox{as}\quad  x\to\pm\infty   
\eeqq
is a bijection onto $\bbr$. As a consequence, the diffusion \eqref{1} is recurrent under $\vth$ with invariant measure 
\beqq\label{inv_ms-allg}
\mu^\vth(dy) \;\;=\;\; \frac{1}{\si^2}\; \frac{1}{s(y)}\; dy  \;\;=\;\; \frac{1}{\si^2}\; \left( \sqrt{1+y^2}\, \right) ^{+\la_1(\vth)} \exp\left( +\sum_{\nu=1}^m \la_{2,\nu}(\vth) F_{2,\nu}(y) \right)
\eeqq
by proposition 9.12~a) in \cite{Ho 14}. \halmos

{\bf Remark 2: } $\Theta_1$ in \eqref{2} is a maximal open parameter interval such that diffusion \eqref{1} is null recurrent:\\
The explicit representation \eqref{klein_s} of $s$, valid with $\vth_1\in\bbr$, shows that $s$ is integrable whenever
 $\vth_1>\frac{\si^2}{2}$ (i.e.\ $\la_1(\vth)>1$): then $S$ in \eqref{harmonic} converges to finite limits $S(\pm\infty)$ as $x\to\pm\infty$, thus $\xi$ is transient by \cite{GS 72}, lemma 3 on p.\ 117.   
\\
In the limiting case $\vth_1=\frac{\si^2}{2}$ ($\la_1(\vth)=1$) to transience, $S$ in \eqref{harmonic}, of logarithmic growth as $x\to\pm\infty$, is a bijection onto $\bbr$. Thus all cases $\vth_1\le\frac{\si^2}{2}$ ($\la_1(\vth)\le 1$) lead to recurrence of $\xi$ with invariant measure $\frac{1}{\si^2}\frac{1}{s(x)}dx$, cf.\ proposition 9.12~a) in \cite{Ho 14}. Ergodicity (positive recurrence) holds if the invariant measure is a finite measure, i.e.\ for $\vth_1<-\frac{\si^2}{2}$ ($\la_1(\vth)<-1$).

{\bf Proof of lemma 2: } See \cite{LS 78} and \cite{JS 87} for local absolute continuity and structure of likelihood ratio processes, a short resum\'e is 6.9--6.12 in \cite{Ho 14}, whence the representation \eqref{log-LR} for log-likelihood ratio processes. Assumption 1 implies that all $f_{2,\nu}$ are bounded functions, thus local $(Q^\vth,\bbf)$-martin\-gales defined by \eqref{score} are square-integrable martingales with angle bracket \eqref{inf_proc} under $Q^\vth$. It remains to prove that $Q^\vth$-almost surely, the process $J(\vth)$  takes values in ${\tt D}^+$, the set of strictly positive definite symmetric matrices in $\bbr^{(1+m)\times(1+m)}$. For $u\in\bbr^{1+m}$ with $|u|=1$, 
\beqq\label{invert}
u^\top J(\vth)_t\; u \;=\;  \frac{1}{\si^2} \int_0^t [u^\top \psi\; \psi^\top u](\eta_s)\,  ds  \;=\; \frac{1}{\si^2} \int_0^t [u^\top \psi]^2(\eta_s)\,  ds \;, 
\quad,\quad 0<t<\infty
\eeqq
where we write for short $\psi:\bbr\to\bbr^{1+m}$ for the row vector with entries $f_1, f_{2,1}, \ldots, f_{2,m}$. Now $Q^\vth$-almost surely, the range $\{\eta_s : 0\le s\le t\}\subset\bbr$ of the diffusion path \eqref{1} has non-empty interior; on this (random, since defined by the path) open set, the functions  $f_1, f_{2,1}, \ldots, f_{2,m}$ are linearly independent by assumption 1. Hence $Q^\vth$-almost surely, the right hand side in \eqref{invert} is strictly positive, first pointwise in $u$ and $t$, second uniformly on $\{|u|=1\}\times [\frac1m,m]$ by continuity, and we let $m\to\infty$. \halmos

{\bf Proof of proposition 1: }  Common determinations $Y$ of the stochastic integrals jointly for all $Q_\vth$ do exist since probability measures $Q^\vth$, $\vth\in\Theta$, are locally equivalent relative to $\bbf$, see lemma 8.2' in \cite{Ho 14}.  Now we exploit almost sure invertibility of the information $J(\vth)_t$ at time $0<t<\infty$ combined with 'inverse bowl shape' of the log-likelihoods. \halmos

{\bf Proof of proposition 2: }  By assumption, the starting point $x_0$ is an interior point of $A$. Thus for $0<t<\infty$,  $\,A\cap\{\eta_s : 0\le s\le t\}$ almost surely contains open balls, and \eqref{invert} remains true with $f{\cdot}1_A$ in place of $f$. \halmos

{\bf Proof for remark 1: } $\check\vth_t$ corresponds to the maximum of the one-dimensional surface 
$$ 
\vth'_1 \quad\lra\quad
(\vth'_1-\vth_1)\, \check S(\vth)_t \;-\; \frac12\, (\vth'_1-\vth_1)^2 \check J_t   \quad,\quad \vth'_1\in\bbr \;. 
$$
The representation of estimation errors under  $\vth\in\Theta$ follows from \eqref{1}, \eqref{score} and the definition of $Y_t$ in proposition 1. By Harris recurrence, the ratio limit theorem holds:  for $g,h\in L^1(\mu^\vth)$ with $\mu^\vth(h)>0$, 
\beqq\label{RLT}
%\mbox{for $g,h\in L^1(\mu^\vth)$ with $\mu^\vth(h)>0$, }\quad 
\frac{ \int_0^t g(\eta_s)\, ds }{ \int_0^t h(\eta_s)\, ds } \quad\lra\quad \frac{ \mu^\vth(g) }{ \mu^\vth(h) } \quad\quad\mbox{ $Q^\vth$-almost surely as $t\to\infty$} \;. 
\eeqq
This is valid for all $\vth\in\Theta$ and for arbitrary choice of a starting point $x_0$ for the process \eqref{1}.   \halmos

%%%%%%%%%%%%%%%%%%%%%%%%%%%%%%%%%%%%%%%%%
%\newpage
\section{Convergence}\label{convergence}
%%%%%%%%%%%%%%%%%%%%%%%%%%%%%%%%%%%%%%%%%

We formulate a theorem on convergence of additive functionals and martingale additive functionals in the null recurrent diffusion 
\eqref{1}. It combines theorem 3.1 from H\"opfner and L\"ocherbach \cite{HL 03} with results due to Khasminskii (\cite{Kh 80}, or theorem 2.2 in \cite{KY 00} and theorem 1.1 in \cite{Kh 01}). The approach is analoguous to \cite{HK 03} or to examples 3.5 and 3.10 in \cite{HL 03}. 

%%%%%%%%%%%%%%%%%
\subsection{Convergence of martingales together with their angle bracket}\label{convergence theorem}
%%%%%%%%%%%%%%%

Introducing further notation for $\vth\in\Theta$, we define 
\beqq\label{psi_p_m}
\Psi^+(\vth)  \;:=\;   \exp\left\{ \,\sum_{\nu=1}^m \la_{2,\nu}\, F_{2,\nu}(+\infty)  \right\} \quad,\quad 
\Psi^-(\vth)  \;:=\;   \exp\left\{ \,\sum_{\nu=1}^m \la_{2,\nu}\, F_{2,\nu}(-\infty)  \right\}    
\eeqq
(cf.\ \eqref{3}, assumption 1 and lemma 1) and introduce a $(0,1)$-valued index  
\beqq\label{alpha}
\al(\vth) \;:=\; \frac12\left( 1 - \frac{2\vth_1}{\si^2} \right) \;=\; \frac12\left( 1 - \la_1(\vth) \right)  
\eeqq
together with weights 
\beqq\label{def_D}
D(\vth) \;:=\; \frac{ (2\si^2)^{1+\al(\vth)}\; \Gamma(\al(\vth)) }{ 2\; \Gamma(1-\al(\vth))  } \;. 
\eeqq
Note that  \eqref{alpha} and \eqref{def_D} depend on $\vth_1\in\Theta_1$ whereas \eqref{psi_p_m} depends on $( \vth_{2,1},\ldots, \vth_{2,m} )\in\bbr^m$. We shall write $\phi:\bbr\to\bbr^k$ for measurable functions whose components $\phi_j$, $1\le j\le k$, are such that  
\beqq\label{def_phi}
\mbox{ $\,|\phi_j|\,$ is locally bounded on $\bbr$, and }\quad 
\int_{\{ |x|>1\}} \phi_j^2(x)\, |x|^{1-\vep}\, dx \;<\; \infty    \quad\mbox{for every $\vep>0$} \;. 
\eeqq
Functions which satisfy \eqref{def_phi} do belong to $L^2(\mu^\vth)$ for every $\vth\in\Theta$, by lemma 1, and we can define 
\beqq\label{def_Sigma}
\Sigma(\vth,\phi) \;:=\; \left( \begin{array}{lll} 
\mu^{\vth}(\,\phi_1^2\,) & \ldots & \mu^{\vth}(\phi_1\phi_k) \\ 
\quad\vdots & \ddots & \quad\vdots \\
\mu^{\vth}(\phi_k\phi_1) &\ldots  & \mu^{\vth}(\,\phi_k^2\,)  
\end{array} \right) \;. 
\eeqq

\vskip0.5cm
{\bf Theorem 1: } Let $\phi$ with components $\phi_j$, $1\le j\le k$, satisfy \eqref{def_phi}. With notations \eqref{psi_p_m}--\eqref{def_D} introduce  sequences of norming constants by 
\beqq\label{norming_seq}
\al_n(\vth) \;:=\; n^{\al(\vth)}\; D(\vth)\,/\, [ \Psi^+(\vth) + \Psi^-(\vth) ] \quad,\quad n\in\bbn \;. 
\eeqq
Then for every $\vth\in\Theta$, we have weak convergence in the Skorohod space $D( \bbr^k , \bbr^{k\times k} )$ of  
\beqq\label{resc_pair_mart_qv}
\left(\; \frac{1}{\sqrt{\al_n(\vth)}} \int_0^{\bullet n} \phi(\eta_s)\, dm^{(\vth)}_s \;\;,\;\; \frac{1}{\al_n(\vth)}  \left\langle \int_0^{\bullet n} \phi(\eta_s)\, dm^{(\vth)}_s\right\rangle \; \right) \quad\quad\mbox{under $Q^\vth$}
\eeqq 
as $\nto$ to 
\beqq\label{limit_pair_mart_qv}
\left(\;  \Sigma^{1/2}\; B\circ V^{\al(\vth)} \;\;,\;\;  \Sigma\;\; V^{\al(\vth)} \;\; \right) \quad,\quad \Sigma = \Sigma(\vth,\phi)
\eeqq
with notation \eqref{def_Sigma}. Here $\,V^\al\,$ is a Mittag-Leffler process of index $0<\al<1$, the process inverse (i.e.\ the process of level crossing times) of a stable increasing process $S^\al$ with  index $0<\al<1$, and $B$ a $k$-dimensional Brownian motion which is independent of $V^\al$; thus $\,B\circ V^{\al}\,$ is Brownian motion subject to independent time change $t\to V^{\al}_t$. \\

The one-dimensional case $k=1$ of theorem 1 corresponds to theorem A in \cite{HK 03} (note that the invariant measure there has a factor $2$  with respect to our $\mu^\vth$ in lemma 1).

{\bf Remark 3: } The one-sided stable process $S^\al$ with index $0<\al<1$ has stationary independent increments with Laplace transform  
$$
E\left( e^{\,-\zeta\, S^\al_t } \right) \;=\; e^{\,-t\;\zeta^\al} \quad,\quad  0<\zeta<\infty\;. 
$$
Its process inverse $V^\al$ has continuous nondecreasing paths with $V^\al_0=0$ and $\lim_{t\uparrow\infty}V^\al_t = \infty$, hence $B\circ V^\al$, $V^\al$ are continuous processes.  
As a consequence (\cite{JS 87}, VI.2.3) of the continuous mapping theorem, theorem 1  implies weak convergence in $\bbr^k{\times}\bbr^{k\times k}$ of \eqref{resc_pair_mart_qv} evaluated at time $t=1$  to \eqref{limit_pair_mart_qv} evaluated at time $t=1$.

{\bf Remark 4: } $\call(V^\al_1)$ is concentrated on $(0,\infty)$. Whereas $\call(V^\al_1)$ admits finite moments of arbitrary order $n\in\bbn$, there is no finite moment of order $-1$: we have 
\beqq\label{neg_moment_ML1}
\int_{(0,\infty)} \call(V^\al_1)(dv)\;\frac1v\;  \;=\; \infty  \quad\mbox{for}\quad 0<\al<1    
\eeqq
which will be of importance below.

%%%%%%%%%%%%%%%%%%%%%%%%%
%\newpage
\subsection{Proof of theorem 1}\label{convergence_proofs}
%%%%%%%%%%%%%%%%%%%%%%%%%%

We formulate a lemma which corresponds to Khasminskii \cite{Kh 80} (see theorem 2.2 in \cite{KY 00} and theorem 1.1 in \cite{Kh 01}, or proposition 9.14 in \cite{Ho 14}).

{\bf Lemma 3: } The path of $\xi$ can be decomposed into iid life cycles  $\xi_{ ]] R_n,R_{n+1}]] }$, $n\ge 1$, defined by 
$$
R_n \;:=\; \inf\{ t>S_n : \xi_t<0 \} \;,\; S_n \;:=\; \inf\{ t>R_{n-1} : \xi>S^{-1}(1) \} \;,\; n\ge 1 \;,\; R_0=0  
$$ 
(up to an initial segment $\xi_{ [[0,R_1]] }$), with $S^{-1}$ the function inverse of $S$ in \eqref{gross_S} which depends on $\vth$. For this decomposition the following holds:  

i) for $\phi$ nonneative and measurable, 
$$
E\left( \int_{R_n}^{R_{n+1}} \phi(\xi_s)\, ds \right) \;\;=\;\;  2\; \mu ^{\vth}(\phi)   \quad\mbox{for all $n\ge 1$} \;;
$$ 

ii) as $t\to\infty$, with $\al=\al(\vth)$ from \eqref{alpha} and $\Psi^{\pm}=\Psi^{\pm}(\vth)$ from \eqref{psi_p_m}:  
$$
P\left( R_{n+1}-R_n > t \right) \quad\sim\quad t^{-\al}\;\; \frac{1}{\Gamma(\al)}\; \left(\frac{1}{2\,\si^2}\right)^\al\, 2\, \left\{ \Psi^+ + \Psi^- \right\} \;. 
$$

{\bf Proof: } 1) The function $S$ in \eqref{harmonic} is harmonic for the Markov generator of the diffusion $\xi$. In a first step we consider $\xi$ transformed by $S$: $\;\wt\xi:=S\circ\xi$ is a diffusion without drift 
$$
d\wt\xi_t \;=\; \wt\si(\wt\xi_t)\, dW_t \quad,\quad 
\wt\si (\wt x) = \si \left( s\circ S^{-1} \right)(\wt x) \quad,\quad \wt x\in\bbr   
$$
(with constant $\si$ in \eqref{1}, and $S^{-1}$ is the function inverse of $S$) whose invariant measure is given by 
\beqq\label{inv_ms_xitilde}
\wt \mu^{\vth} (d\wt x)\; =\; \frac{1}{ [\wt\si(\wt x) ]^2}\; d\wt x  \quad\mbox{on}\quad (\bbr,\calb(\bbr))   
\eeqq
(\cite{GS 72}, \cite{Kh 80}, or proof of proposition 9.12 in \cite{Ho 14}). Calculating from \eqref{klein_s} and \eqref{gross_S} the function $s\circ S^{-1}$ and supressing the dependence on $\vth$ we get
$$
\left( s\circ S^{-1} \right)(\wt x) \;\;\sim\;\;  |\wt x |^{\frac{-\la_1}{1-\la_1}}\; (1-\la_1)^{\frac{-\la_1}{1-\la_1}}\;  (\Psi^{\pm})^{\frac{-1}{1-\la_1}} \quad\mbox{as}\quad \wt x \to\pm\infty
$$
with notation from \eqref{3} and \eqref{psi_p_m}, thus with notation \eqref{alpha}
$$
\frac{1}{ [\wt\si(\wt x) ]^2}  \;\;\sim\;\;  \frac{1}{\si^2}\;  \left|\wt x \right|^{ \la_1 / \al }\; (1-\la_1)^{ \la_1 / \al }\; (\Psi^{\pm})^{1/\al} \quad\mbox{as}\quad \wt x \to\pm\infty  \;. 
$$
This is the Lebesgue density of the invariant measure $\wt \mu^{\vth}$ in \eqref{inv_ms_xitilde} for the process $\,\wt\xi=S\circ\xi\,$.  

2) We comment on the invariant density for $\wt\xi$. Multiplying the asymptotic constant by $2$ we define 
\beqq\label{A_pm}
A^\pm \;:=\;  \frac{2}{\si^2} (1-\la_1)^{ \la_1 / \al }\; (\Psi^{\pm})^{1/\al}   
\eeqq
and have 
\beqq\label{inv_dens_xitilde}
\frac{2}{ [\wt\si(\wt x) ]^2}  \;\;\sim\;\;  \left|\wt x \right|^{ \la_1 / \al }\;  A^{\pm} \quad\mbox{as}\quad \wt x \to\pm\infty  \;.   
\eeqq
This is the norming used by Khasminskii for the invariant measure. By \eqref{3} and \eqref{alpha}, the quotient 
$$
\beta \;:=\; \la_1 / \al = 2 \frac{\la_1}{1-\la_1}  
$$
takes values in $(-1,\infty)$ and is such that
$$
\frac{1}{\beta+2} \;=\; \frac{ \frac12(1-\la_1) }{ \la_1 + (1-\la_1) }  \;=\; \frac12(1-\la_1) \;=\; \al \;. 
$$

3) Now we apply the results of \cite{Kh 80} quoted in the beginning of section \ref{convergence_proofs}. The stopping times   
$$
R_n \;:=\; \inf\{ t>S_n : \wt\xi_t<0 \} \;,\; S_n \;:=\; \inf\{ t>R_{n-1} : \wt\xi>1 \} \;,\; n\ge 1 \;,\; R_0=0  
$$ 
define iid life cycles $\wt\xi_{ ]] R_n,R_{n+1}]] }$ in the path of $\wt\xi$ for $n\ge 1$, with the following two properties: first,    
\beqq\label{prop_1}
E\left( \int_{R_n}^{R_{n+1}} \wt\phi(\wt\xi_s)\, ds \right) \;=\;  2\; \wt \mu ^{\vth}(\wt\phi)   \quad,\quad n\ge 1 
\eeqq
for measurable functions $\wt\phi:\bbr\to[0,\infty)$ and $\wt\mu^\vth$ as in \eqref{inv_ms_xitilde}; second,  
\beqq\label{prop_2}
P\left( R_{n+1}-R_n > t \right) \;\;\sim\;\; t^{-\al}\; \frac{\;\left[ \al\right]^{2\al}}{\Gamma(1+\al)}  \left\{ [A^+]^\al + [A^-]^\al \right\} 
\quad\mbox{as}\quad t\to\infty \;. 
\eeqq
4)  The iid life cycles for $\wt\xi= (S(\xi_t))_{t\ge 0}$ in step 3) are iid life cycles for $\xi=(\xi_t)_{t\ge 0}$ 
$$
R_n \;:=\; \inf\{ t>S_n : \xi_t<0 \} \;,\; S_n \;:=\; \inf\{ t>R_{n-1} : \xi>S^{-1}(1) \} \;,\; n\ge 1 \;,\; R_0=0  
$$
for which asymptotics  \eqref{prop_2} remain unchanged, and change of variables transforms \eqref{prop_1} into 
$$
E\left( \int_{R_n}^{R_{n+1}} \phi(\xi_s)\, ds \right) \;=\;  2\; \mu ^{\vth}(\phi)   \quad\mbox{for all $n\ge 1$} \;.
$$
Using step 2), the asymptotic constant on the right hand side of \eqref{prop_2} equals 
$$
\frac{\;\left[ \al\right]^{2\al}}{\Gamma(1+\al)}  \left\{ [A^+]^\al + [A^-]^\al \right\} \;\;=\;\;  
\frac{1}{\Gamma(\al)}\, \left(\frac{1}{2\si^2}\right)^\al\,  2\; \left\{ \Psi^+ + \Psi^- \right\} 
$$
and lemma 3 is proved by step 3); recall that both $\al=\al(\vth)$ and $\la_1=\la_1(\vth)$ depend on $\vth\in\Theta$. \halmos\\

{\bf Proof of theorem 1: } 1) In a first step, consider convergence of martingales in case $k=1$. \\Fix $\vth\in\Theta$ and  apply  theorem 3.1~c) of \cite{HL 03} to one-dimensional measurable functions $\phi$ which belong to $L^2(\mu^\vth)$: to the sequence $(R_n)_n$ of renewal times considered in lemma 3 correspond norming sequences 
$$
\wt a_n \;=\; \frac{1}{ \Gamma(1-\al) \; P( R_2-R_1 > n ) } \;\;=\;\; n^{\al}\; \frac{\Gamma(\al)}{\Gamma(1-\al)}\, \left(2 \si^2 \right)^\al  \frac{1}{ 2\, \{ \Psi^+ + \Psi^-  \} }
$$
%($\wt a_n$, $\al$ and $\la_1$ depending on $\vth\in\Theta$) 
such that 
$$
 \frac{1}{\sqrt{\wt \al_n}} \int_0^{\bullet n} \phi(\eta_s)\, dm^{(\vth)}_s 
$$
converges weakly as $\nto$ in the Skorohod space $D$ of one-dimensional cadlag functions to 
$$
\sqrt{\wt c\,}\; B\circ V^{\al(\vth)} 
$$
where $\wt c=\wt c(\vth)$ is given by 
\beqq\label{def_ctilde}
\wt c  \;\;:=\;\;  E_\vth\left(\, \int_{R_1}^{R_2}  \phi^2(\eta_s)\;  d \langle m^{(\vth)} \rangle_s \,\right)  \;=\; 2\si^2\; \mu^\vth(\phi^2)  \;; 
\eeqq
here we use \eqref{prop_1} together with  $d \langle m^{(\vth)} \rangle_s = \si^2\, ds$. 
Taking into account this factor $2\si^2$ arising on the right hand side of \eqref{def_ctilde} we arrive at the norming sequence $(\al_n)_n$ of theorem 1,  defined by \eqref{norming_seq} and \eqref{def_D}, and have weak convergence  under $Q^\vth$  
\beqq\label{1d_mart_convergence}
 \frac{1}{\sqrt{\al_n}} \int_0^{\bullet n} \phi(\eta_s)\, dm^{(\vth)}_s  \quad\lra\quad \sqrt{c\,}\; B\circ V^{\al(\vth)} 
\quad\mbox{(weakly in $D$ as $\nto$)}  
\eeqq 
with $\,c := \mu^\vth(\phi^2)\,$ as asserted in theorem 1. This proves convergence of martingales in case $k=1$. 

2) The result proved in step 1) can be extended to case $k>1$ and to weak convergence of martingales together with their angle bracket: apply   corollaries 3.2 and 3.3 in \cite{HL 03}.\halmos\\

{\bf Proof for remark 3: }  We refer to 2.5-2.8 in \cite{HL 03} and the references quoted there. \halmos

{\bf Proof for remark 4: } It is well known that $E( [S^\al_1]^p )$ is finite for $0<p<\al$ and infinite for $p=\al$ (e.g.\ \cite{Ho 14}, 6.18'). Since $V^\al$ is the process inverse to $S^\al$, the relation 
$$
\call(\,V^\al_1 \,) \;=\; \call(\, \left[\frac{1}{S^\al_1}\right]^\al ) \quad,\quad 0<\al<1    
$$
holds and proves that $E( \frac{1}{[V^\al_1]} )$ equals $\infty$.\halmos

%%%%%%%%%%%%%%%%%%%%%%%%%%%%%%%%%%%%%%%%%
%\newpage
\section{LAMN and optimal estimator sequences}\label{LAMN}
%%%%%%%%%%%%%%%%%%%%%%%%%%%%%%%%%%%%%%%%%

For local asymptotic mixed normality see Jeganathan \cite{Je 82}, Davies \cite{Da 85}, and LeCam and Yang \cite{LY 90}. See also sections 5.1, 6.1 and 7 of \cite{Ho 14}. 

%%%%%%%%%%%%%%%%%
\subsection{Limit distribution for the estimators of section \ref{assumptions} }
%%%%%%%%%%%%%%%%%

Under $\vth\in\Theta$ consider the norming sequence $( \al_n(\vth) )_n$ defined by  \eqref{norming_seq}. 
In notations from \eqref{score} and \eqref{inf_proc}, proposition 1 combined with theorem 1 yields convergence in law as $\nto$ of rescaled ML errors
\beqq\label{MLE_conv}
\call\left( \sqrt{ \al_n(\vth)\, } \left( \wh\vth _n - \vth \right) \mid Q^\vth \,\right)  \quad\lra\quad \call\left(\, \Lambda(\vth)^{-1/2}\, B(V^{\al(\vth)}_1)/V^{\al(\vth)}_1  \,\right) 
\eeqq  
for all $\vth\in\Theta$, where $B$ is $(1{+}m)$-dimensional Brownian motion independent of the Mittag-Leffler variable $V^{\al(\vth)}_1$, and where $\Lambda(\vth)$ is given by
\beqq\label{limit_Lambda}
\Lambda(\vth) \;=\; \frac{1}{\si^4} \left( \begin{array}{lll} 
\mu^{\vth}(\,f_1^2\,) & \ldots & \mu^{\vth}(f_1 f_{2,\nu'}) \\ 
\quad\vdots & \ddots & \quad\vdots \\
\mu^{\vth}( f_{2,\nu} f_1 ) &\ldots  & \mu^{\vth}(\,f_{2,\nu}f_{2,\nu'}\,)  
\end{array} \right)   
\eeqq
with $\nu$, $\nu'$ ranging from $1$ to $m$.  If we fix as in proposition 2 a compact interval $A$ such that $x_0\in {\rm int}(A)$ and replace %the functions  
$\,f_1$, $f_{2,1},\ldots,f_{2,m}\,$ by $\,f_1{\cdot}1_A$, $f_{2,1}{\cdot}1_A,\ldots,f_{2,m}{\cdot}1_A\,$, 
we get  convergence in law as $\nto$  
\beqq\label{est_conv}
\call\left( \sqrt{ \al_n(\vth)\, } \left( \wt\vth(A) _n - \vth \right) \mid Q^\vth \,\right)  \quad\lra\quad \call\left(\, \Lambda(\vth,A)^{-1/2}\, B(V^{\al(\vth)}_1)/V^{\al(\vth)}_1  \,\right) 
\eeqq
where $\Lambda(\vth,A)$ is given by 
\beqq\label{limit_Lambda_A}
\Lambda(\vth,A) \;=\; \frac{1}{\si^4} \left( \begin{array}{lll} 
\mu^{\vth}(\,f_1^2\,{\cdot}1_A\,) & \ldots & \mu^{\vth}(f_1 f_{2,\nu'}\, {\cdot}1_A) \\ 
\quad\vdots & \ddots & \quad\vdots \\
\mu^{\vth}( f_{2,\nu} f_1\, {\cdot}1_A ) &\ldots  & \mu^{\vth}( f_{2,\nu}f_{2,\nu'}\,{\cdot} 1_A )  
\end{array} \right)   
\eeqq
for $1\le \nu,\nu'\le m$. Note that the deterministic matrices in \eqref{limit_Lambda} and  \eqref{limit_Lambda_A} are invertible for all $\vth\in\Theta$. The argument is similiar to \eqref{invert}: with $\psi:\bbr\to\bbr^{1+m}$ defined by its components $\,f_1$, $f_{2,1},\ldots, f_{2,m}\,$  (or by $\,f_1{\cdot}1_A$, $f_{2,1}{\cdot}1_A,\ldots,f_{2,m}{\cdot}1_A\,$ where $A$ contains open balls), the strictly positive Lebesgue density of $\mu^\vth$ in lemma~1 together with assumption 1~iii) shows that  
$$
u^\top \Lambda(\vth)\; u \;=\;  \int [u^\top \psi\; \psi^\top u]\,  d\mu^{\vth}   \;=\; \int [u^\top \psi]^2\, d\mu^{\vth} \;\;>\;\; 0 
$$
for $u\in\bbr^{1+m}$ with $|u|=1$. This is half-ordering of symmetric nonnegative nonnegative definite matrices, and the argument above also shows  $\,\Lambda(\vth) > \Lambda(\vth,A)\,$ in this sense. Thus the limit law in \eqref{est_conv} is necessarily more spread out than the limit law in \eqref{MLE_conv}. On the other hand,  independence of $B$ and $V^\al$ allows to write       
\beqq\label{limit_structure_ML}
\call \left( B(V^{\al(\vth)}_1)/V^{\al(\vth)}_1  \right) \quad=\quad   \int_{(0,\infty)} \call\left( V^{\al(\vth)}_1 \right)(dv)\; \caln\left( 0 \,,\, \frac1v\, I \right) 
\eeqq   
with $I$ the identity matrix in $\bbr^{(1+m)\times(1+m)}$ and $0\in\bbr^{1+m}$. \eqref{limit_structure_ML} implies that the limit laws in  \eqref{MLE_conv} or \eqref{est_conv} never admit finite variances: by \eqref{neg_moment_ML1} in remark 4, we have  
$$
\int_{(0,\infty)} \call\left( V^{\al(\vth)}_1 \right)(dv)\; \frac1v \;\;=\;\; \infty    
$$
for all $\vth\in\Theta$. Given the structure of the limit laws in \eqref{est_conv} or \eqref{MLE_conv}, classical theory of estimation based on comparison of estimators through variances (or finite moments of higher order) of limit distributions for rescaled estimation errors is of no use for the statistical model \eqref{1}  under consideration. 
%for the estimators considered in section 1.1. 
%in view of limit distibutions as in \eqref{est_conv} or \eqref{MLE_conv}.  

%%%%%%%%%%%%%%%%%%%%%%%%%%%%%5
\subsection{LAMN and comparison of estimators locally at $\vth\in\Theta$}\label{comp_est}
%%%%%%%%%%%%%%%%%%%%%%%%%%%%%%%

For $\vth\in\Theta$ we reconsider the norming sequence $( \al_n(\vth) )_n$ of \eqref{norming_seq} and call    
\beqq\label{local_scale}
\delta_n(\vth) \;:=\;  \frac{1}{\sqrt{n^{\al(\vth)}\,}}  \quad,\quad 
\al(\vth)=\frac12(1-\frac{2 \vth_1}{\si^2}) \quad\mbox{ranging over}\;\; (0,1)
\eeqq 
local scale at $\vth$. The following combines lemma 2 with theorem 1: if we write 
$$
{\bf S}_n(\vth) \;:=\;   \delta_n(\vth)\; S(\vth)_n  %\;=\;  \frac{1}{\sqrt{\al_n(\vth)\,}}\; S(\vth)_n 
\quad,\quad 
{\bf J}_n(\vth) \;:=\;   \delta_n^2(\vth)\; J_n %\;=\; \frac{1}{\al_n(\vth)}\; J_n 
$$
for score and information rescaled with \eqref{local_scale} and call  
$$
\cale(\vth,n) \;:=\; \left\{ Q^{\vth + \delta_n(\vth) h} \,:\;  h \in \Theta_{\vth,n} \right\} 
$$
local experiment at $\vth$, parametrized by $h$ as local parameter for which  
$$
\Theta_{\vth,n} \;:=\; \left\{ h\in\bbr^{1+m} : \vth + \delta_n(\vth) h \in \Theta \right\} 
$$
increases to $\bbr^{1+m}$ as $\nto$, then log-likelihoods in $\cale(\vth,n)$   
\beqq\label{limit_logLRs}
{\bf \Lambda}_{n,\vth}^{ h / 0 } \;:=\;   \Lambda_n^{ ( \vth + \delta_n(\vth) h ) / \vth }   \;=\; h^\top {\bf S}_n(\vth) \;+\;  h^\top {\bf J}_n(\vth) \; h \quad,\quad  h \in \Theta_{\vth,n}
\eeqq
are quadratic in the local parameter $h$. Theorem 1 with norming sequence \eqref{norming_seq} establishes weak convergence of pairs as $\nto$ 
\beqq\label{def_bfS_bfJ}
\call\left( \left( {\bf S}_n(\vth) \,,\, {\bf J}_n(\vth) \right) \mid Q^\vth \right) \quad\lra\quad 
\left( {\bf S}(\vth) \,,\, {\bf J}(\vth) \right) \;:=\; 
\left(\,  {\bf \Sigma}(\vth)^{1/2}\; B\circ V^{\al(\vth)} \;,\;  {\bf \Sigma}(\vth)\;\; V^{\al(\vth)} \, \right)  
\eeqq
where  %${\bf \Sigma}(\vth)$ in \eqref{def_bfS_bfJ} denotes the deterministic matrix
\beqq\label{def_bfSigma}
{\bf \Sigma}(\vth) \;\;:=\;\; \frac{ D(\vth) }{ \Psi^+(\vth) + \Psi^-(\vth) } \;\cdot\; \Lambda(\vth) 
\eeqq
combines  scaling factors from \eqref{psi_p_m}--\eqref{def_D}  with the deterministic and invertible $(1{+}m){\times}(1{+}m)$-matrix 
$$
\Lambda(\vth) \;=\; \frac{1}{\si^4} \left( \begin{array}{lll} 
\mu^{\vth}(\,f_1^2\,) & \ldots & \mu^{\vth}(f_1 f_{2,\nu'}) \\ 
\quad\vdots & \ddots & \quad\vdots \\
\mu^{\vth}( f_{2,\nu} f_1 ) &\ldots  & \mu^{\vth}(\,f_{2,\nu}f_{2,\nu'}\,)  
\end{array} \right)   
$$
(with $\nu$, $\nu'$ ranging from $1$ to $m$) which was defined in \eqref{limit_Lambda} above.

By proposition 1, $\,J_n\,$ is taking values in ${\tt D}^+$, $\,Q^\vth$-almost surely for all $\vth\in\Theta$: thus the same holds for ${\bf J}_n(\vth)$ defined above. As a consequence of log-likelihoods \eqref{limit_logLRs} and convergence \eqref{def_bfS_bfJ} in the local model $\cale(\vth,n)$ at $\vth$, local asymptotic mixed normality (LAMN) holds at every point $\vth\in\Theta$, with central sequence given by 
\beqq\label{central_seq}
{\bf Z}_n(\vth)  \;:=\; 1_{\{ {\bf J}_n(\vth) \in {\tt D}^+ \}}\,  {\bf J}_n^{-1}(\vth)\, {\bf S}_n(\vth) \;=\; 1_{\{ J_n \in {\tt D}^+ \}}\, \delta_n^{-1}(\vth) J_n^{-1}\, S(\vth)_n 
\quad,\quad n\in\bbn \;.   
\eeqq
The limit experiment generated by the pair $( {\bf S}(\vth) , {\bf J}(\vth) )$  admits     
\beqq\label{central_stat}
{\bf Z}(\vth) \;:=\;  {\bf J}^{-1}(\vth)\, {\bf S}(\vth) \;=\;  {\bf \Sigma}(\vth)^{-1/2}\; B( V^{\al(\vth)}_1 ) /  V^{\al(\vth)}_1  
\eeqq 
as central statistic. At every point $\vth\in\Theta$, we thus dispose of the following local asymptotic minimax bound (see \cite{LY 90} theorem 1 in section 5.6, or \cite{Ho 14} theorem 7.12): \\

{\bf Theorem 2: } (LeCam) Consider any loss function $\ell:\bbr^{1+m}\to[0,\infty)$ which is continuous, bounded and subconvex (i.e.\ level sets $\{ x : \ell(x)\le c \}$ are convex and symmetric in the sense $\ell(x)=\ell(-x)$, $x\in\bbr^{1+m}$). Fix any point $\vth\in\Theta$. 

a) Consider any sequence of $\calf_n$-measurable estimators $T_n$ for the unknown parameter such that 
\beqq\label{tightness_cond}
\call\left( \delta_n^{-1}(\vth) \left(T_n-\vth\right) \mid Q^\vth \right) \;,\; n\ge 1\;,\quad\mbox{is tight in $\bbr^{1+m}$ as $\nto$}      
\eeqq
where $(\delta_n(\vth))_n$ is local scale \eqref{local_scale}. Then 
\beqq\label{lamm_bound}
\lim_{c\to\infty}\;\liminf_{\nto}\; \sup_{|h|\le c}\; E_{Q^{\vth+\delta_n(\vth)h}}\left(\; \ell\left(\; \delta_n^{-1}(\vth)(T_n-(\vth+\delta_n(\vth)h))\;\right) \;\right) \quad\ge\quad E\left(\;\ell({\bf Z(\vth)})\;\right)   
\eeqq
where ${\bf Z}(\vth)$ is the central statistic \eqref{central_stat} in the limit experiment.  

b) If rescaled estimation errors at $\vth$ can be coupled to the central sequence \eqref{central_seq} at $\vth$ 
\beqq\label{coupling_cond}
\delta_n^{-1}(\vth) \left(T_n-\vth\right) \quad=\quad {\bf Z}_n(\vth) \;+\; o_{Q^\vth}(1)\quad\mbox{as $\nto$}   \;, 
\eeqq
then the estimator sequence  $(T_n)_n$ achieves    
$$
\lim_{\nto}\; \sup_{|h|\le c}\; E_{Q^{\vth+\delta_n(\vth)h}}\left(\; \ell\left(\; \delta_n^{-1}(\vth)(T_n-(\vth+\delta_n(\vth)h))\;\right) \;\right) \quad=\quad E\left(\;\ell({\bf Z(\vth)})\;\right) 
$$ 
for arbitrary constants $c<\infty$, and in particular attains the local asymptotic minimax bound \eqref{lamm_bound}. \\

Combining the representation \eqref{estim_errors_1} of rescaled ML errors in proposition 1 to \eqref{coupling_cond}, we see that the ML sequence $(\wh\vth_n)_n$ satisfies the coupling condition. Hence theorem 2 states that at every point $\vth$ in the parameter space $\Theta$, the sequence $(\wh\vth_n)_n$ minimizes errors uniformly over shrinking neighbourhoods of $\vth$ whose radius is proportional to local scale $\delta_n(\vth)$ as $\nto$. Thus we dispose of an explicit and tractable estimator sequence which is asymptotically optimal in the sense of the local asymptotic minimax theorem.

The estimators $\,\wt \theta (A)_n\,$ considered in proposition 2 do not satisfy the coupling condition \eqref{coupling_cond}, but converge at rate $\sqrt{n^{\al(\vth)}}$ as $\nto$. Under some conditions in addition to  LAMN (see Davies \cite{Da 85}, or \cite{Ho 14} section 7.4) one step modification due to LeCam allows to transform preliminary estimators with the right speed of convergence into modified ones which do satisfy the  coupling condition \eqref{coupling_cond}, and thus are optimal in the sense of the convolution theorem. In our case however, the situation is much simpler: the log-likelihood ratios in \eqref{log-LR} are exactly quadratic in the parameter, the information \eqref{inf_proc} is free from the parameter, and mimicking the relation which according to \eqref{1} relates the score at $\vth$ to $Y_n$ of proposition 1
$$
S_n(\vth) \;=\; Y_n \;-\; J_n\; \vth
$$
an elementary definition of a 'score with estimated parameter' 
$$
S_n( T_n ) \;:=\; Y_n \;-\; J_n\; T_n
$$
is possible whenever $T_n$ is $\calf_n$-measurable and $\bbr^{1+m}$-valued. Then $Q^\vth$-almost surely 
$$
T_n  \;+\;  1_{\{ J_n \in {\tt D}^+ \}}\,  J_n^{-1}\, S_n( T_n ) \;\;=\;\; T_n  \;+\; 1_{\{ J_n \in {\tt D}^+ \}}\left\{  J_n^{-1}\, Y_n  \;-\; T_n \right\}  \;\;=\;\; \wh\vth_n \;+\; 1_{\{ J_n \notin {\tt D}^+ \}}\, T_n
$$
is the ML estimator, for every $\vth\in\Theta$, thus 'one-step correction' works in our model \eqref{1} in its most elementary form.

%%%%%%%%%%%%%%%%%%%%%%%%%%%%%%%
%\newpage
\bibliographystyle{plainnat}

\vskip2.0cm
{\tt
Reinhard H\"opfner, Institut f\"ur Mathematik, Universit\"at Mainz \\
Staudingerweg 9, 55099 Mainz\\
hoepfner@mathematik.uni-mainz.de\\

Carina Zeller, Institut f\"ur Mathematik, Universit\"at Mainz\\
Staudingerweg 9, 55099 Mainz\\
czeller@students.uni-mainz.de}

%%%%%%%%%%%%%%%%%%%%%%%%%%%%%%%%%%%%
\end{document}